# A THEORETICAL COMPARISON OF THE DATA AUGMENTATION, MARGINAL AUGMENTATION AND PX-DA ALGORITHMS

By James P. Hobert[1] and Dobrin Marchev[2]

*University of Florida and Baruch College, CUNY*

The data augmentation (DA) algorithm is a widely used Markov chain Monte Carlo (MCMC) algorithm that is based on a Markov transition density of the form $p(x|x') = \int_Y f_{X|Y}(x|y) f_{Y|X}(y|x') \, dy$, where $f_{X|Y}$ and $f_{Y|X}$ are conditional densities. The PX-DA and marginal augmentation algorithms of Liu and Wu [*J. Amer. Statist. Assoc.* **94** (1999) 1264–1274] and Meng and van Dyk [*Biometrika* **86** (1999) 301–320] are alternatives to DA that often converge much faster and are only slightly more computationally demanding. The transition densities of these alternative algorithms can be written in the form $p_R(x|x') = \int_Y \int_Y f_{X|Y}(x|y') R(y, dy') f_{Y|X}(y|x') \, dy$, where $R$ is a Markov transition function on $Y$. We prove that when $R$ satisfies certain conditions, the MCMC algorithm driven by $p_R$ is at least as good as that driven by $p$ in terms of performance in the central limit theorem and in the operator norm sense. These results are brought to bear on a theoretical comparison of the DA, PX-DA and marginal augmentation algorithms. Our focus is on situations where the group structure exploited by Liu and Wu is available. We show that the PX-DA algorithm based on Haar measure is at least as good as any PX-DA algorithm constructed using a proper prior on the group.

**1. Introduction.**

1.1. *Background.* In statistical problems where there is a need to explore an intractable density, $f_X(x)$, there is sometimes available a joint density

Received April 2006; revised February 2007.
[1]Supported by NSF Grants DMS-00-72827 and DMS-05-03648.
[2]Supported by NSF Grant DMS-00-72827 and PSC-CUNY Grant 60075-35 36.

$f(x,y)$, on $\mathsf{X} \times \mathsf{Y}$ say, such that $\int_{\mathsf{Y}} f(x,y)\,dy = f_X(x)$ and such that simulating from the conditional densities, $f_{X|Y}(x|y)$ and $f_{Y|X}(y|x)$, is straightforward. In such situations, one can apply the data augmentation (DA) algorithm (Tanner and Wong [22]), which is a Markov chain Monte Carlo (MCMC) algorithm based on the Markov transition density (Mtd) given by

$$(1) \qquad p(x|x') = \int_{\mathsf{Y}} f_{X|Y}(x|y) f_{Y|X}(y|x')\,dy.$$

It is well known and easy to show that $p(x|x')$ is reversible with respect to $f_X$, which implies that $f_X$ is an invariant density. Like its cousin, the EM algorithm, the DA algorithm is considered a useful algorithm that sometimes suffers from slow convergence.

The PX-DA algorithm (Liu and Wu [13]) and the closely related marginal augmentation (MA) algorithm (Meng and van Dyk [14]) are alternatives to DA that often converge much faster and are only slightly more computationally demanding. The basic idea is to use $f(x,y)$ to create an entire family of joint densities that all have $f_X$ as the $x$ marginal. Each member of this family can be used to form a DA algorithm and the hope is that some of the resulting algorithms will be significantly better than the original. To be specific, consider a class of functions $t_g: \mathsf{Y} \to \mathsf{Y}$ for $g \in G$ such that, for each fixed $g$, $t_g(y)$ is one-to-one and differentiable in $y$. Suppose further that $r(g)$ is a probability density on $G$ and define another probability density $\tilde{f}: \mathsf{X} \times \mathsf{Y} \times G \to [0, \infty)$ as $\tilde{f}(x,y,g) = f(x, t_g(y))|J_g(y)|r(g)$, where $J_g(z)$ is the Jacobian of the transformation $z = t_g^{-1}(y)$. Let $\tilde{f}(x,y) = \int_G \tilde{f}(x,y,g)\,dg$ and note that $\int_{\mathsf{Y}} \tilde{f}(x,y)\,dy = f_X(x)$. The PX-DA algorithm (which is the same as the MA algorithm in this situation) is simply the alternative DA algorithm based on the Mtd given by

$$(2) \qquad p_r(x|x') = \int_{\mathsf{Y}} \tilde{f}_{X|Y}(x|y) \tilde{f}_{Y|X}(y|x')\,dy.$$

By varying $r(\cdot)$, we can create the family of joint densities mentioned above. Liu and Wu [13], Meng and van Dyk [14] and van Dyk and Meng [23] (hereafter, L&W, M&vD and vD&M) have provided many examples where this strategy leads to major improvements over standard DA algorithms.

Straightforward sampling from $\tilde{f}_{X|Y}$ and $\tilde{f}_{Y|X}$, which is necessary if the PX-DA algorithm is to be useful in practice, is made possible by exploiting the relationship between these conditionals and the joint density $\tilde{f}(x,y,g)$. First, consider sampling from $\tilde{f}_{Y|X}$ and note that

$$\tilde{f}_{Y|X}(y|x) = \int_G f_{Y|X}(t_g(y)|x)|J_g(y)|r(g)\,dg.$$

Consequently, we can draw from $\tilde{f}_{Y|X}$ by drawing $y'$ and $g$ independently from $f_{Y|X}(y'|x)$ and $r(g)$, respectively, and setting $y = t_g^{-1}(y')$. Now let $f_Y(y) = \int_X f(x,y)\,dx$ and let $w(g;y)$ denote the density proportional to $r(g)|J_g(y)|f_Y(t_g(y))$. We can draw from $\tilde{f}_{X|Y}$ by drawing $g$ from $w(g;y)$ and then $x \sim f_{X|Y}(x|t_g(y))$. Putting all of this together, as in [13], Scheme 1.1, a single iteration of the PX-DA algorithm ($x' \to x$) can be accomplished by performing the following three steps:

1. Draw $y \sim f_{Y|X}(y|x')$.
2. Draw $g \sim r(\cdot)$, draw $g'$ from $w(g'; t_g^{-1}(y))$ and set $y' = t_{g'}(t_g^{-1}(y))$.
3. Draw $x \sim f_{X|Y}(x|y')$.

Note that the first and third steps are exactly the same as the two steps of the DA algorithm. Given that $w(g;y)$ contains the term $f_Y(t_g(y))$ and that direct sampling from $f_Y$ is infeasible (otherwise MCMC would be unnecessary), one might expect that sampling from $w(g;y)$ would be difficult. However, as the examples in [13, 14, 23] illustrate, when $g$ has lower dimension than $y$, sampling from $w(g;y)$ can be completely straightforward, adding very little to the overall computational burden.

We use Albert and Chib's [1] DA algorithm for Bayesian probit regression as a running example. Let $V_1, V_2, \ldots, V_n$ denote independent random variables with $V_i \mid \beta \sim \text{Bernoulli}(\Phi(z_i^T \beta))$ where $z_i$ is a $p \times 1$ vector of known covariates associated with $V_i$, $\beta$ is a $p \times 1$ vector of unknown regression coefficients and $\Phi(\cdot)$ is the standard normal distribution function. A flat prior on $\beta$ leads to an (intractable) posterior density given by

$$\pi(\beta|v) = \frac{1}{m(v)} \prod_{i=1}^n [\Phi(z_i^T \beta)]^{v_i} [1 - \Phi(z_i^T \beta)]^{1-v_i},$$

where $m(v)$ is the marginal mass function. Let $\mathbb{R}^+ = (0, \infty)$, $\mathbb{R}^- = (-\infty, 0]$ and consider the function

$$\pi(\beta, y|v) = \frac{1}{m(v)} \left[ \prod_{i=1}^n \{I_{\mathbb{R}^+}(y_i) I_{\{1\}}(v_i) + I_{\mathbb{R}^-}(y_i) I_{\{0\}}(v_i)\} \phi(y_i; z_i^T \beta, 1) \right],$$

where $y = (y_1, y_2, \ldots, y_n)^T \in \mathbb{R}^n$, $I_A(\cdot)$ is the indicator of the set $A$ and $\phi(x; \mu, \sigma^2)$ denotes the $N(\mu, \sigma^2)$ density function evaluated at the point $x$. Straightforward calculations show that $\pi(\beta, y|v)$ is a joint density in $(\beta, y)$ whose $\beta$ marginal is the target, $\pi(\beta|y)$. Moreover, $\pi(\beta|y, v)$ is a multivariate normal density and $\pi(y|\beta, v)$ is a product of $n$ truncated univariate normal densities. Albert and Chib's algorithm alternates between these two conditionals. L&W developed a PX-DA algorithm for this problem by taking $t_g(y) = gy$ and $G = (0, \infty)$. This yields $w(g;y) \propto r(g)g^n \exp\{-g^2 y^T My/2\} I_{\mathbb{R}^+}(g)$, where $M$ is a known $n \times n$ matrix. Drawing from the multivariate density $\pi(y|v)$ does not appear straightforward, but sampling from the univariate density $w(g;y)$ is easy as long as $r(g)$ has a simple form. Indeed, L&W take $r(g) \propto g^{a-1} e^{-bg^2} I_{\mathbb{R}^+}(g)$ where $a, b > 0$, which allows one to sample from $w(g;y)$ by drawing a gamma variate and taking the square root.

1.2. *A general class of alternatives to DA.* Step 2 of the PX-DA algorithm involves making the transition $y \to y'$ and can therefore be interpreted as simulating one step of a Markov chain on $\mathsf{Y}$. In fact, Theorem 1 in [13] shows that $f_Y$ is an invariant density for this chain. Thus, the Mtd of the PX-DA algorithm is a special case of the general Mtd given by

$$(3) \qquad p_R(x|x') = \int_\mathsf{Y} \int_\mathsf{Y} f_{X|Y}(x|y') R(y, dy') f_{Y|X}(y|x') \, dy,$$

where $R(y, dy')$ is any Markov transition function (Mtf) on $\mathsf{Y}$ that has $f_Y$ as an invariant density. Routine calculations show that $f_X$ is invariant for $p_R$ and that, if $R$ is reversible with respect to $f_Y$, then $p_R$ is reversible with respect to $f_X$. In this paper, we perform the first general study of (3). The main results provide conditions under which the Markov chain driven by (3) is better than the corresponding DA algorithm. To be specific, we show that if $p_R$ is reversible with respect to $f_X$, then $p_R$ is at least as good as $p$ in the efficiency ordering of Mira and Geyer [16], which concerns performance in the central limit theorem (CLT). (For a cleaner exposition, we henceforth write "better than" instead of the more accurate "at least as good as.") We also show that if $p_R$ is itself a DA algorithm; that is, if there exists a joint density $f^*(x, y)$ such that $\int_\mathsf{Y} f^*(x, y) \, dy = f_X(x)$ and such that $p_R$ can be reexpressed as

$$p_R(x|x') = \int_\mathsf{Y} f^*_{X|Y}(x|y) f^*_{Y|X}(y|x') \, dy,$$

then $p_R$ is better than $p$ in the operator norm sense (Liu, Wong and Kong [11]).

Because the Mtds of the DA, MA and PX-DA algorithms can all be written in the form (3), our general results concerning (3) can be brought to bear on a theoretical comparison of these algorithms. This yields both new results and generalizations of known results from [13, 14] and [23]. Furthermore, our proofs of the generalizations are simpler and require fewer regularity conditions than the original proofs. It is our hope that the results herein will promote theoretical and methodological development of improved DA algorithms.

Here is a simple example of the application of our results concerning (3). The PX-DA algorithm is, by definition, a DA algorithm and as such is reversible with respect to $f_X$. Hence, the results described above are applicable and imply that every PX-DA algorithm is better than the DA algorithm in the efficiency ordering and in the operator norm sense. The efficiency ordering result is new, but the operator norm result is known—see Theorem 2 in [13] and Theorem 1 in [14]. Note that we say "every PX-DA algorithm." This is because the result holds no matter what (proper) density $r(g)$ is used to construct the PX-DA algorithm.

1.3. *Adapting to an improper $r(g)$: Liu and Wu's group structure.* L&W, M&vD and vD&M all argued that the PX-DA algorithm should perform better as the density $r(g)$ becomes more "diffuse" or "spread out," and they provided empirical evidence supporting this claim. It is clearly impossible to implement the PX-DA algorithm in the limiting case where $r$ is improper. However, L&W and M&vD found (what appear to be) different ways of utilizing an improper $r(g)$ to construct an algorithm that achieves the limiting convergence rate. L&W developed their results by exploiting a certain group structure that may be present in the problem. M&vD, on the other hand, constructed a nonpositive recurrent Markov chain on $\mathsf{X} \times G$ having stationary density $f_X(x)r(g)$ and provided conditions under which the $x$ component of that chain is itself a Markov chain with invariant density $f_X$. We focus on L&W's approach and show that, when the group structure exists, L&W's algorithm is exactly the same as M&vD's algorithm (under a particular improper working prior). This is the first formal comparison of the two limiting algorithms. We now briefly describe L&W's group structure and limiting algorithm.

Suppose that $G$ is a topological group; that is, a group such that the functions $(g_1, g_2) \mapsto g_1 g_2$ and $g \mapsto g^{-1}$ are both continuous. Let $e$ denote the group's identity element. (An example is the *multiplicative group*, $\mathbb{R}^+$, where group composition is defined as multiplication, the identity element is $e = 1$ and $g^{-1} = 1/g$.) Suppose further that $t_e(y) = y$ for all $y \in \mathsf{Y}$ and that $t_{g_1 g_2}(y) = t_{g_1}(t_{g_2}(y))$ for all $g_1, g_2 \in G$ and all $y \in \mathsf{Y}$. Assume that $G$ is a unimodular group and let $\nu(dg)$ denote Haar measure on $G$. One iteration of L&W's limiting algorithm, which we call the Haar PX-DA algorithm, consists of the following three steps:

1. Draw $y \sim f_{Y|X}(y|x')$.
2. Draw $g$ from the density (with respect to $\nu$) proportional to $|J_g(y)| f_Y(t_g(y))$ and set $y' = t_g(y)$.
3. Draw $x \sim f_{X|Y}(x|y')$.

Note that the Haar PX-DA algorithm actually requires less computation than the PX-DA algorithm. Indeed, Step 2 involves only a single draw from a distribution on $G$, while the middle step of the PX-DA algorithm requires two such draws. The Mtd associated with this algorithm has $f_X$ as an invariant density (see L&W) and is, in fact, another special case of (3). (Note that the invariance of $f_X$ is not obvious in this case because, unlike PX-DA, the Haar PX-DA algorithm is not defined as an alternative DA algorithm.) L&W proved that the Haar PX-DA algorithm is better in the operator norm sense than every PX-DA algorithm.

Consider again the probit regression example. The multiplicative group, $G = \mathbb{R}^+$, is unimodular with Haar measure given by $\nu(dg) = dg/g$ where $dg$ denotes Lebesgue measure. Furthermore, the transformation $t_g(y) = gy$ satisfies the compatibility conditions described above, so the Haar PX-DA algorithm is applicable. As shown in [13], the middle step entails drawing $g$ from a density proportional to $g^{n-1}\exp\{-g^2 y^T M y\}I_{\mathbb{R}^+}(g)$. Both L&W and vD&M provide strong empirical evidence that this algorithm can converge much faster than Albert and Chib's [1] DA algorithm.

1.4. *Comparing general versions of PX-DA and Haar PX-DA.* We develop generalizations of the PX-DA and Haar PX-DA algorithms in a setting where $X$, $Y$ and $G$ are abstract spaces (not necessarily Euclidean) and the group $G$ is not required to be unimodular (Haar measure is replaced by left-Haar measure). This is accomplished in two steps. First, the group structure is used to build Mtfs, $Q_r(y, dy')$ and $Q(y, dy')$, that are reversible with respect to $f_Y$ and that behave like general versions of the middle steps of the PX-DA and Haar PX-DA algorithms. Then Mtds for the generalized versions of PX-DA and Haar PX-DA are formed by using $Q_r$ and $Q$ in place of $R$ in (a generalized version of) (3). Because L&W did not use the term "Haar PX-DA," it is important to bear in mind throughout this paper that what we call the "general Haar PX-DA algorithm" is, in fact, a generalization of L&W's limiting PX-DA algorithm.

A comparison of the resulting generalized algorithms is facilitated by a representation of Haar PX-DA as an improvement of PX-DA. More specifically, we show that there exists a joint density $\tilde{f}(x, y)$, whose $x$ marginal is $f_X$, such that the Mtd of the general PX-DA algorithm can be written as

$$
(4) \quad \begin{aligned} p_r(x|x') &= \int_\mathsf{Y} \int_\mathsf{Y} f_{X|Y}(x|y) Q_r(y, dy') f_{Y|X}(y|x') \mu_y(dy) \\ &= \int_\mathsf{Y} \tilde{f}_{X|Y}(x|y) \tilde{f}_{Y|X}(y|x') \mu_y(dy), \end{aligned}
$$

where $\mu_x(dx)$ and $\mu_y(dy)$ are the analogues of $dx$ and $dy$ that will be defined in Section 3. (This, of course, implies that PX-DA is better than DA.) We then show that the Mtd of the general Haar PX-DA algorithm can be written as

$$
\begin{aligned} p^*(x|x') &= \int_\mathsf{Y} \int_\mathsf{Y} f_{X|Y}(x|y) Q(y, dy') f_{Y|X}(y|x') \mu_y(dy) \\ &= \int_\mathsf{Y} \int_\mathsf{Y} \tilde{f}_{X|Y}(x|y) \tilde{Q}(y, dy') \tilde{f}_{Y|X}(y|x') \mu_y(dy), \end{aligned}
$$

where $\tilde{f}$ is as in (4) and $\tilde{Q}(y, dy')$ is reversible with respect to $\int_\mathsf{X} \tilde{f}(x, y) \mu_x(dx) = \tilde{f}_Y(y)$; that is, $p^*$ is an improvement of $p_r$. It is also shown that $p^*(x|x')$ is itself a DA algorithm. Therefore, our results concerning (3) imply that $p^*(x|x')$ is better than every version of $p_r(x|x')$ in the efficiency ordering and in the operator norm sense. As before, the efficiency ordering result is new, but a special case of the operator norm result was established in [13] (see Section 5 for details).

The remainder of the paper is laid out as follows. In Section 2, we set notation and review some results from general state space Markov chain theory. Our study of (3) commences in Section 3. In Section 4, we describe two different methods of using a group action to construct a Mtf with a prespecified stationary distribution. Finally, our general versions of the PX-DA and Haar PX-DA algorithms are introduced and studied in Section 5.

**2. Markov chain background.** As in Meyn and Tweedie ([15], Chapter 3) let $P(x, dy)$ be a Mtf on a set $\mathsf{X}$ equipped with a countably generated $\sigma$-algebra $\mathcal{B}(\mathsf{X})$. Suppose that $\pi$ is an invariant probability measure; that is, $\pi(A) = \int_\mathsf{X} P(x, A)\pi(dx)$ for all measurable $A$. Denote the Markov chain defined by $P(x, dy)$ as $\Phi = \{\Phi_n\}_{n=0}^\infty$, where the distribution of $\Phi_0$ will be stated explicitly when needed. As usual, let $L^2(\pi)$ be the vector space of real-valued, measurable functions on $\mathsf{X}$ that are square-integrable with respect to $\pi$, and let $L_0^2(\pi)$ be the subspace of mean zero functions; that is, functions satisfying $\int_\mathsf{X} f(x)\pi(dx) = 0$. Define inner product on this space by $\langle f, g \rangle = \int_\mathsf{X} f(x)g(x)\pi(dx)$. The corresponding norm is given by $\|f\| = \sqrt{\langle f, f \rangle}$. The Mtf $P(x, dy)$ defines an operator, $P$, that acts on $f \in L_0^2(\pi)$ through

$$(Pf)(x) = \int_\mathsf{X} f(y)P(x, dy) = \mathrm{E}[f(\Phi_{n+1})|\Phi_n = x].$$

Note that $\langle Pf, f \rangle = \mathrm{Cov}(f(\Phi_0), f(\Phi_1))$ when $\Phi_0 \sim \pi$. The chain $\Phi$ (or, equivalently, the Mtf $P$) is said to be reversible with respect to $\pi$ if for all bounded functions $f, g \in L_0^2(\pi)$, $\langle Pf, g \rangle = \langle f, Pg \rangle$. The norm of the operator $P$ is defined as

$$\|P\| = \sup_{f \in L_0^2(\pi), f \neq 0} \frac{\|Pf\|}{\|f\|} = \sup_{f \in L_0^2(\pi), \|f\|=1} \|Pf\|.$$

A straightforward application of Jensen's inequality shows that $\|P\| \leq 1$.

Now assume that $\int_\mathsf{X} |h(x)|\pi(dx) < \infty$ and that MCMC will be used to estimate the intractable expectation $\pi h := \int_\mathsf{X} h(x)\pi(dx)$. If $\Phi$ is irreducible, aperiodic and Harris recurrent (see Meyn and Tweedie [15] for definitions), then the ergodic average $\overline{h}_n = n^{-1}\sum_{i=0}^{n-1} h(\Phi_i)$ converges almost surely to $\pi h$ no matter what the distribution of $\Phi_0$. This justifies the use of $\overline{h}_n$ as an estimator of $\pi h$. There are several different methods available for calculating the standard error of this estimator (see, e.g., Geyer [5], Hobert, Jones, Presnell and Rosenthal [7] and Jones, Haran, Caffo and Neath [8]) and all are based on the assumption that there is a CLT for $\overline{h}_n$; that is, that there exists a $\sigma^2 \in (0, \infty)$ such that, as $n \to \infty$, $\sqrt{n}(\overline{h}_n - \pi h) \xrightarrow{d} \mathrm{N}(0, \sigma^2)$. The asymptotic variance, $\sigma^2$, depends on both the function $h$ and the Mtf $P$ (but not on the distribution of $\Phi_0$) so we write it as $v(h, P)$. If the CLT fails to hold, then we simply write $v(h, P) = \infty$.

Unfortunately, even if $h \in L^2(\pi)$, irreducibility, aperiodicity and Harris recurrence (henceforth "the usual regularity conditions") are not enough to guarantee that $v(h, P) < \infty$. The chain is called *geometrically ergodic* if there exist $M : \mathsf{X} \to [0, \infty)$ and $\rho \in [0, 1)$ such that $\|P^n(x, \cdot) - \pi(\cdot)\|_{\mathrm{TV}} \leq M(x) \rho^n$ for all $x \in \mathsf{X}$ and all $n = 1, 2, 3, \ldots$, where $\|\cdot\|_{\mathrm{TV}}$ denotes total variation norm. If $\Phi$ is geometrically ergodic and reversible with respect to $\pi$, then $v(h, P) < \infty$ for every $h \in L^2(\pi)$ (Roberts and Rosenthal [17]). Many popular Monte Carlo Markov chains have been shown to be geometrically ergodic. See, for example, Jones and Hobert [9] and Roberts and Rosenthal [18], and the references therein.

Now suppose that we wish to estimate $\pi h$ and we have available two different Mtf's, $P$ and $Q$, with invariant probability measure $\pi$ such that $v(h, P)$ and $v(h, Q)$ are both finite. If $P$ and $Q$ are similar in terms of simulation effort, then we would clearly prefer the more efficient chain; that is, the chain with the smaller asymptotic variance. Moreover, if $v(h, P) \leq v(h, Q)$ for all $h$, then we would prefer $P$ over $Q$ regardless of the function $h$. This discussion motivates the following definitions from Mira and Geyer [16].

DEFINITION 1. If $P$ and $Q$ are two Mtf's with invariant probability measure $\pi$ that both satisfy the usual regularity conditions, then $P$ is better than $Q$ in the efficiency ordering, written $P \succeq_E Q$, if $v(h, P) \leq v(h, Q)$ for every $h \in L^2(\pi)$.

DEFINITION 2. If $P$ and $Q$ are two Mtf's with invariant probability measure $\pi$, then $P$ dominates $Q$ in the covariance ordering, written $P \succeq_1 Q$, if $\langle Ph, h \rangle \leq \langle Qh, h \rangle$ for every $h \in L_0^2(\pi)$.

The following result provides a characterization of the efficiency ordering for reversible chains as well as a practical method of proving that $P \succeq_E Q$.

THEOREM 1 (Mira and Geyer [16]). *Let $P$ and $Q$ be two Mtf's that are reversible with respect to the probability measure $\pi$ and that satisfy the usual regularity conditions. Then $P \succeq_E Q$ if and only if $P \succeq_1 Q$.*

It is important to note that $\succeq_E$ provides only a partial ordering; that is, it can happen that neither $P \succeq_E Q$ nor $Q \succeq_E P$ holds. In such a case, neither chain is better than the other and the choice between $P$ and $Q$ will depend on the particular function to be estimated.

Monte Carlo Markov chains can also be compared via their operator norms. Indeed, the quantity $\|P\|$ is closely related to the convergence rate of the corresponding Markov chain. For instance, if $P$ is reversible with respect to $\pi$ and satisfies the usual regularity conditions, then $P$ is geometrically ergodic if and only if $\|P\| < 1$ (Roberts and Rosenthal [17] and Roberts and Tweedie [20]). Furthermore, results in Liu, Wong and Kong [12] show that the smaller the norm, the faster the chain converges. Examples of the use of this criterion for comparing Monte Carlo Markov chains can be found in [11, 13, 14].

It is important to keep in mind that neither $\|P\| \leq \|Q\|$ nor $P \succeq_E Q$ guarantees that $P$ is a good Monte Carlo Markov chain. Indeed, even if $\|P\| \leq \|Q\|$, it may be the case that both $P$ and $Q$ are bad chains (with norm 1) and neither should be used. Similarly, $P \succeq_E Q$ tells us nothing about the existence of CLTs for $P$. However, if $P$ is also known to be geometrically ergodic, then we could rule out $Q$ *and* be content to use $P$ to explore the target distribution. The results described above imply that if $P$ and $Q$ are both reversible and $\|P\| \leq \|Q\|$, then geometric ergodicity of $Q$ implies that of $P$. (See Roberts and Rosenthal [19] for some related results.) This result can be extremely useful in practice because the better chain ($P$ in this case) is typically more complex and hence harder to analyze. This idea is exploited in Roy and Hobert [21], who prove that the Haar PX-DA algorithm for the probit model (discussed in Section 1) is geometric by showing that the simpler DA algorithm of Albert and Chib is geometric.

**3. Improving upon the DA algorithm.** In this section, we study Mtds of the form (3). Assume that $\mathsf{X}$ and $\mathsf{Y}$ are locally compact, separable metric spaces equipped with their Borel $\sigma$-algebras. Assume further that $\mu_x$ and $\mu_y$ are $\sigma$-finite measures on $\mathsf{X}$ and $\mathsf{Y}$, respectively, and that $f(x,y)$ is a probability density on $\mathsf{X} \times \mathsf{Y}$ with respect to $\mu_x \times \mu_y$. As usual, let $f_X$, $f_Y$, $f_{X|Y}$ and $f_{Y|X}$ denote the marginal and conditional densities. In this context, the DA algorithm has Mtd (with respect to $\mu_x$) given by

$$(5) \qquad p(x|x') = \int_{\mathsf{Y}} f_{X|Y}(x|y) f_{Y|X}(y|x') \mu_y(dy).$$

The analogue of (3) is

$$(6) \qquad p_R(x|x') = \int_{\mathsf{Y}} \int_{\mathsf{Y}} f_{X|Y}(x|y') R(y, dy') f_{Y|X}(y|x') \mu_y(dy),$$

where $R(y, dy')$ is any Mtf on $\mathsf{Y}$ that has $f_Y$ as an invariant density. Again, straightforward calculations reveal that $f_X$ is an invariant density for $p_R$ and that reversibility of $R$ with respect to $f_Y$ implies reversibility of $p_R$ with respect to $f_X$. Varying the Mtf $R(y, dy')$ produces a family of Markov chains having $f_X$ as invariant density, and (as we explain later) the DA algorithm is one of the family members. In some cases, $p_R$ is itself a DA algorithm; that is, there exists a probability density $f^*(x,y)$ on $\mathsf{X} \times \mathsf{Y}$ with respect to $\mu_x \times \mu_y$ such that $\int_\mathsf{Y} f^*(x,y)\mu_y(dy) = f_X(x)$ and such that $p_R$ can be reexpressed as

$$p_R(x|x') = \int_\mathsf{Y} f^*_{X|Y}(x|y) f^*_{Y|X}(y|x') \mu_y(dy).$$

Clearly, if $p_R$ is a DA algorithm, then it is reversible with respect to $f_X$. We now state a known result about DA that will be used to prove the main result in this section.

THEOREM 2 (Amit [2] and Liu, Wong and Kong [11]). *Let $P$ denote the operator on $L_0^2(f_X)$ associated with $p(x|x')$. Let $(X,Y) \sim f(x,y)$ and $h \in L_0^2(f_X)$. Then $\langle Ph, h \rangle = \mathrm{Var}[\mathrm{E}(h(X)|Y)]$ and $\|P\| = \gamma^2(X,Y)$, where $\gamma(X,Y)$ is the maximal correlation between $X$ and $Y$.*

The next result allows us to compare two different versions of (6).

THEOREM 3. *Suppose that $R$ and $S$ are two Mtf's on $\mathsf{Y}$ that have $f_Y$ as invariant density and assume that $R \succeq_1 S$. Let $p_R$ and $p_S$ denote the corresponding versions of (6) and denote the associated operators as $P_R$ and $P_S$. Assume that $p_R$ and $p_S$ satisfy the usual regularity conditions. If $p_R$ and $p_S$ are both reversible with respect to $f_X$, then $p_R \succeq_E p_S$. If, in addition, $p_R$ and $p_S$ are both DA algorithms, then $\|P_R\| \leq \|P_S\|$.*

PROOF. Let $\Phi^* = \{\Phi_n^*\}_{n=0}^\infty$ and $\tilde{\Phi} = \{\tilde{\Phi}_n\}_{n=0}^\infty$ denote stationary versions of the chains driven by $p_R$ and $p_S$, respectively. Fix $h \in L_0^2(f_X)$ and define $h^*(y) = \int_\mathsf{X} h(x) f_{X|Y}(x|y) \mu_x(dx)$. It is easy to see that $h^* \in L_0^2(f_Y)$. Now

$$\begin{aligned}\langle P_R h, h \rangle &= \int_\mathsf{X} \int_\mathsf{X} h(x') h(x) p_R(x|x') f_X(x') \mu_x(dx') \mu_x(dx) \\ &= \int_\mathsf{X} \int_\mathsf{X} \int_\mathsf{Y} \int_\mathsf{Y} h(x') h(x) f_{X|Y}(x|y') R(y, dy') f_{Y|X}(y|x') f_X(x') \\ &\quad \times \mu_y(dy) \mu_x(dx') \mu_x(dx) \\ &= \int_\mathsf{Y} \int_\mathsf{Y} \left[\int_\mathsf{X} h(x) f_{X|Y}(x|y') \mu_x(dx)\right] \left[\int_\mathsf{X} h(x') f_{X|Y}(x'|y) \mu_x(dx')\right] \\ &\quad \times R(y, dy') f_Y(y) \mu_y(dy) \\ &= \int_\mathsf{Y} \int_\mathsf{Y} h^*(y) h^*(y') R(y, dy') f_Y(y) \mu_y(dy) \\ &\leq \int_\mathsf{Y} \int_\mathsf{Y} h^*(y) h^*(y') S(y, dy') f_Y(y) \mu_y(dy) = \langle P_S h, h \rangle,\end{aligned}$$

where the inequality follows from the fact that $R \succeq_1 S$. It then follows from Theorem 1 that $p_R \succeq_E p_S$. Now let $f^*(x,y)$ and $\tilde{f}(x,y)$ denote the densities that allow us to express $p_R$ and $p_S$ as DA algorithms. In conjunction with the results above, Theorem 2 implies that

$$\text{Var}[\text{E}(h(X^*)|Y^*)] = \langle P_R h, h \rangle \leq \langle P_S h, h \rangle = \text{Var}[\text{E}(h(\tilde{X})|\tilde{Y})],$$

where $(X^*, Y^*) \sim f^*(x,y)$ and $(\tilde{X}, \tilde{Y}) \sim \tilde{f}(x,y)$. Now, since $X^* \stackrel{d}{=} \tilde{X}$, we have $\{g : \text{Var}\, g(X^*) = 1\} = \{g : \text{Var}\, g(\tilde{X}) = 1\}$. Suppose that $\text{Var}\, g(X^*) = 1$ and put $\mu_g = \int_{\mathsf{X}} g(x) f_X(x)\, dx$. Then

$$\begin{aligned}(7)\quad \text{Var}[\text{E}(g(X^*)|Y^*)] &= \text{Var}\{\text{E}[(g(X^*) - \mu_g)|Y^*]\} \\ &\leq \text{Var}\{\text{E}[(g(\tilde{X}) - \mu_g)|\tilde{Y}]\} = \text{Var}[\text{E}(g(\tilde{X})|\tilde{Y})].\end{aligned}$$

But it is well know that for random elements $U$ and $V$,

$$\gamma^2(U, V) = \sup_{\{h\, :\, \text{Var}\, h(U) = 1\}} \text{Var}[\text{E}(h(U)|V)].$$

It follows that $\|P_R\| = \gamma^2(X^*, Y^*) \leq \gamma^2(\tilde{X}, \tilde{Y}) = \|P_S\|$. □

Theorem 3 actually allows us to compare the DA algorithm with the algorithm based on (6). Indeed, the Mtd (5) can be viewed as a special case of (6) where $R(y, dy')$ is taken to be the trivial Mtf that is a point mass at $y$. This trivial Mtf is obviously dominated in the covariance ordering by any nontrivial $R$. We conclude that, if $p_R$ can be expressed as a DA algorithm, then it is better than the original DA algorithm both in terms of efficiency and operator norm. We state this as a corollary.

COROLLARY 1. *Suppose that $R$ is a Mtf on $\mathsf{Y}$ that has $f_Y$ as invariant density. Let $p_R$ be as in (6) and denote the associated operator by $P_R$. Assume that $p$ and $p_R$ satisfy the usual regularity conditions. If $p_R$ is reversible with respect to $f_X$, then $p_R \succeq_E p$. If, in addition, $p_R$ is a DA algorithm, then $\|P_R\| \leq \|P\|$.*

In order to apply Corollary 1, we must establish that $p_R$ is reversible and possibly that $p_R$ is a DA algorithm. We know that reversibility of $R$ implies that of $p_R$. The next result shows that there is also a simple condition on $R$ that implies that $p_R$ is a DA algorithm.

PROPOSITION 1. *Let $R$ be a Mtf on $\mathsf{Y}$ that has $f_Y$ as invariant density and let $p_R$ be as in (6). If there exists a Mtf $R^{1/2}(y, dy')$ that is reversible with respect to $f_Y$ and is such that $R(y, dy') = \int_{\mathsf{Y}} R^{1/2}(w, dy') R^{1/2}(y, dw)$, then $p_R$ is a DA algorithm with respect to $f^*(x,y) = f_Y(y) \int_{\mathsf{Y}} f_{X|Y}(x|y') R^{1/2}(y, dy')$.*

PROOF. First, it is easy to see (without using reversibility) that $\int_{\mathsf{Y}} f^*(x,y) \times \mu_y(dy) = f_X(x)$. Now

$$\int_{\mathsf{Y}} f^*_{X|Y}(x|y) f^*_{Y|X}(y|x') \mu_y(dy)$$

$$= \int_{\mathsf{Y}} \left[ \frac{f^*(x,y)}{\int_{\mathsf{X}} f^*(x,y)\mu_x(dx)} \right] \left[ \frac{f^*(x',y)}{\int_{\mathsf{Y}} f^*(x',y)\mu_y(dy)} \right] \mu_y(dy)$$

$$= \int_{\mathsf{Y}} \left[ \int_{\mathsf{Y}} f_{X|Y}(x|y') R^{1/2}(y,dy') \right]$$

$$\quad \times \left[ \frac{f_Y(y)}{f_X(x')} \int_{\mathsf{Y}} f_{X|Y}(x'|y'') R^{1/2}(y,dy'') \right] \mu_y(dy)$$

$$= \int_{\mathsf{Y}} \int_{\mathsf{Y}} \int_{\mathsf{Y}} f_{X|Y}(x|y') R^{1/2}(y'',dy')$$

$$\quad \times \frac{1}{f_X(x')} f_{X|Y}(x'|y) f_Y(y) R^{1/2}(y,dy'') \mu_y(dy)$$

$$= \int_{\mathsf{Y}} \int_{\mathsf{Y}} f_{X|Y}(x|y') f_{Y|X}(y|x') \left[ \int_{\mathsf{Y}} R^{1/2}(y,dy'') R^{1/2}(y'',dy') \right] \mu_y(dy'')$$

$$= p_R(x|x'). \qquad \square$$

Two situations where the hypotheses of Proposition 1 are clearly satisfied are (i) if $R$ is reversible with respect to $f_Y$ and *idempotent* in the sense that $R(y,dy') = \int_{\mathsf{Y}} R(w,dy') R(y,dw)$, and (ii) if $R$ is defined to be the Mtf corresponding to two consecutive steps of a chain on $\mathsf{Y}$ that is reversible with respect to $f_Y$.

**4. Using group actions to construct Markov transition functions.** We now use the group structure on $G$ to build two Mtf's, $Q_r(y,dy')$ and $Q(y,dy')$, that behave like general versions of the middle steps (Step 2) of the PX-DA and Haar PX-DA algorithms described in Section 1.

4.1. *The group structure.* Let $\mathsf{Y}$ and $f_Y$ be as defined in the previous section and assume now that $G$ is another locally compact, separable metric space that is also a topological group. Suppose that the group $G$ acts topologically on the left of $\mathsf{Y}$; that is, there is a continuous function $F:G \times \mathsf{Y} \to \mathsf{Y}$ such that $F(e,y) = y$ for all $y \in \mathsf{Y}$ and $F(g_1 g_2, y) = F(g_1, F(g_2, y))$ for all $g_1, g_2 \in G$ and all $y \in \mathsf{Y}$. [Note that $F(g,y)$ is playing the role of $t_g(y)$ from Section 1.] As is typically done, we will abbreviate $F(g,y)$ with $gy$ so, for example, the second condition is written $(g_1 g_2)y = g_1(g_2 y)$.

As in Eaton [3], we use the term *multiplier* to describe a continuous homomorphism of $G$ into the multiplicative group $\mathbb{R}^+$; that is, a function $\chi:G\to\mathbb{R}^+$ is a multiplier if $\chi$ is continuous and $\chi(g_1g_2)=\chi(g_1)\chi(g_2)$ for all $g_1,g_2\in G$. Clearly, if $\chi$ is a multiplier, then $\chi(e)=1$ and $\chi(g^{-1})=1/\chi(g)$. The measure $\mu_y$ is called *relatively* (*left*) *invariant* with multiplier $\chi$ if

$$\chi(g)\int_{\mathsf{Y}} h(gy)\mu_y(dy)=\int_{\mathsf{Y}} h(y)\mu_y(dy),$$

for all $g\in G$ and all integrable functions $h:\mathsf{Y}\to\mathbb{R}$. As an example, consider the PX-DA algorithm for the probit model that was discussed in Section 1. In that case, the group acts on the left of $\mathsf{Y}=\mathbb{R}^n$ through scalar multiplication, $(g,y)\mapsto gy$, and $\mu_y$, which is Lebesgue measure on $\mathbb{R}^n$, is easily seen to be relatively invariant with multiplier $\chi(g)=g^n$.

While all of the examples considered in [13, 14] and [23] satisfy the assumptions of the previous two paragraphs, this level of generality is not quite enough. In order to ensure that our results subsume those of L&W, we assume that there exists a function $j:G\times\mathsf{Y}\to\mathbb{R}^+$ such that:

1. $j(g^{-1},y)=\frac{1}{j(g,y)}$ $\forall g\in G,\ y\in\mathsf{Y}$,
2. $j(g_1g_2,y)=j(g_1,g_2y)j(g_2,y)$ $\forall g_1,g_2\in G,\ y\in\mathsf{Y}$, and
3. For all $g\in G$ and all integrable functions $h:\mathsf{Y}\to\mathbb{R}$,

$$\int_{\mathsf{Y}} h(gy)j(g,y)\mu_y(dy)=\int_{\mathsf{Y}} h(y)\mu_y(dy). \quad (8)$$

Note that when $\mu_y$ is relatively invariant, we can simply take $j(g,y)$ to be $\chi(g)$. Now suppose (as in [13]) that $\mathsf{Y}\subset\mathbb{R}^n$, $\mu_y$ is Lebesgue measure on $\mathsf{Y}$, and for each fixed $g\in G$, $F(g,\cdot):\mathsf{Y}\to\mathsf{Y}$ is differentiable. Then if we take $j(g,y)$ to be the Jacobian of the transformation $y\mapsto F(g,y)$, the three properties listed above follow straightforwardly from calculus.

4.2. *A transition based on a probability measure on $G$.* We now build a Mtf, $Q_r$, that is a generalized version of Step 2 of the PX-DA algorithm. Let $r$ be a probability measure on $G$. Define

$$m_r(y)=\int_G f_Y(gy)j(g,y)r(dg)$$

and assume that $m_r(y)>0$ for all $y\in\mathsf{Y}$. Define $N=\{y\in\mathsf{Y}:m_r(y)=\infty\}$ and let $\overline{\mathsf{Y}}=\mathsf{Y}\setminus N$. Note that $\int_{\mathsf{Y}} m_r(y)\mu_y(dy)=1$, which implies that $\mu_y(N)=0$. Assume that $gy\in\overline{\mathsf{Y}}$ for all $y\in\overline{\mathsf{Y}}$ and all $g\in G$. A simple calculation shows that, for fixed $y\in\overline{\mathsf{Y}}$,

$$f_Y(g'g^{-1}y)j(g',g^{-1}y)/m_r(g^{-1}y)$$

is a probability density function on $G\times G$ with respect to $r\times r$. Let $Q_r$ be an operator on $L_0^2(f_Y)$ defined as

$$(Q_r h)(y)=\int_G\int_G \frac{h(g'g^{-1}y)f_Y(g'g^{-1}y)j(g',g^{-1}y)}{m_r(g^{-1}y)}r(dg)r(dg')$$

when $y \in \overline{\mathsf{Y}}$ and $(Q_r h)(y) = \int_{\mathsf{Y}} h(y) f_Y(y) \mu_y(dy)$ when $y \in N$. This is the operator corresponding to a Markov chain on $\mathsf{Y}$ that evolves as follows. If the current state, $y$, is in $\overline{\mathsf{Y}}$, then the distribution of the next state is that of $g'g^{-1}y$ where $(g, g')$ is a random element from the density $f_Y(g'g^{-1}y)j(g', g^{-1}y)/m_r(g^{-1}y)$, and if $y \in N$, then the next state is from $f_Y$. Denote the corresponding Mtf on $\mathsf{Y}$ as $Q_r(y, dy')$. We now establish that $f_Y$ is an invariant density for $Q_r$ by showing that $Q_r$ is reversible with respect to $f_Y$.

PROPOSITION 2. *Suppose $r$ is a probability measure on $G$ such that $m_r(y) > 0$ for all $y \in \mathsf{Y}$ and such that $gy \in \overline{\mathsf{Y}}$ for all $y \in \overline{\mathsf{Y}}$ and all $g \in G$. Then the Mtf $Q_r$ is reversible with respect to $f_Y$.*

PROOF. We prove the result in the case where $\mu_y$ is relatively invariant and leave the extension to the general case to the reader. Let $h_1, h_2 \in L_0^2(f_Y)$ be bounded. We will show that $\langle Q_r h_1, h_2 \rangle = \langle h_1, Q_r h_2 \rangle$. Indeed,

$$\langle h_1, Q_r h_2 \rangle = \int_G \int_G \int_{\overline{\mathsf{Y}}} \frac{h_1(y) h_2(g'g^{-1}y) f_Y(y) f_Y(g'g^{-1}y) \chi(g')}{m_r(g^{-1}y)}$$
(9)
$$\times \mu_y(dy) r(dg) r(dg').$$

Now, since $gg'^{-1}g'g^{-1} = e$, the inner integral in (9) can be expressed as

$$\int_{\overline{\mathsf{Y}}} \frac{h_1(gg'^{-1}g'g^{-1}y) h_2(g'g^{-1}y) f_Y(gg'^{-1}g'g^{-1}y) f_Y(g'g^{-1}y) \chi(g'g^{-1}) \chi(g)}{m_r(g'^{-1}g'g^{-1}y)}$$
$$\times \mu_y(dy),$$

which, using the relative invariance of $\mu_y$, becomes

$$\int_{\overline{\mathsf{Y}}} \frac{h_1(gg'^{-1}y) h_2(y) f_Y(gg'^{-1}y) f_Y(y) \chi(g)}{m_r(g'^{-1}y)} \mu_y(dy).$$

Thus, (9) can be written as

$$\int_G \int_G \int_{\overline{\mathsf{Y}}} \frac{h_2(y) h_1(g'g^{-1}y) f_Y(y) f_Y(g'g^{-1}y) \chi(g')}{m_r(g^{-1}y)} \mu_y(dy) r(dg) r(dg')$$
$$= \langle Q_r h_1, h_2 \rangle. \qquad \square$$

EXAMPLE 1. Let $\mathsf{Y} = \mathbb{R}$ and take $\mu_y$ to be Lebesgue measure. Let $f_Y(y) = \frac{1}{2} e^{-|y|}$ and take $G$ to be the multiplicative group on $\mathbb{R}^+$. If the group action is defined as multiplication, then $\mu_y$ is relatively invariant with multiplier $\chi(g) = g$. [We always use $\chi(g)$ instead of $j(g, y)$ when $\mu_y$ is relatively invariant.] If we take $r(dg)$ to be a probability measure with density $e^{-g}$ on the positive half-line, then $m_r(y) = (1 + |y|)^{-2} \in (0, \infty)$ for all $y \in \mathsf{Y}$. [For an example where $m_r$ is not finite everywhere, use $(1 + g)^{-2}$ in place of $e^{-g}$.] A

A COMPARISON OF DATA AUGMENTATION ALGORITHMS    15...

simple calculation shows that the distribution of the random element $(g, g')$ used to make the transitions under $Q_r$ can be described as follows. First, $g \sim \text{Exp}(1)$ and, conditional on $g$, $g'$ has density (with respect to Lebesgue measure on $\mathbb{R}^+$) given by

$$\frac{f_Y(g'g^{-1}y)\chi(g')e^{-g'}}{m_r(g^{-1}y)} \propto g'e^{-g'(1+|y|/g)}.$$

Hence, $g'|g \sim \text{Gamma}(2, 1+|y|/g)$ and it follows that $g'g^{-1} \stackrel{d}{=} v$, where $v$ is a random variable on $\mathbb{R}^+$ with density given by

$$f(v) = v\exp\{-v|y|\}\left[\frac{2}{(v+1)^3} + \frac{2|y|}{(v+1)^2} + \frac{|y|^2}{(v+1)}\right].$$

Consequently, for measurable $A \subset \mathbb{R}$, $Q_r(y, A) = \int_A q_r(y'|y)\mu_y(dy')$ where

$$q_r(y'|y) = e^{-|y'|}\frac{|y'||y|}{|y'+y|}\left[\frac{2}{|y'+y|^2} + \frac{2}{|y'+y|} + 1\right]$$
$$\times [I_{\mathbb{R}^+}(y)I_{\mathbb{R}^+}(y') + I_{\mathbb{R}^-}(y)I_{\mathbb{R}^-}(y')].$$

Clearly, $q_r(y'|y)f_Y(y)$ is a symmetric function of $(y', y)$ so the Mtf $Q_r$ is reversible with respect to $f_Y$ as it must be according to Proposition 2. Note that the chain is not irreducible. For example, if it is started with $y_0 > 0$, then it will never visit the negative half-line.

4.3. *A transition based on left-Haar measure on $G$.* In this section, we build on results in Liu and Sabatti [10] to construct a Mtf, $Q$, that is a generalized version of Step 2 of the Haar PX-DA algorithm. We begin by describing left-Haar measure and some of its properties. Under the assumptions of Section 4.1 there exists a left-Haar measure, $\nu_l$, on $G$, which is a nontrivial measure satisfying

(10) $$\int_G h(\tilde{g}g)\nu_l(dg) = \int_G h(g)\nu_l(dg)$$

for all $\tilde{g} \in G$ and all integrable functions $h: G \to \mathbb{R}$. This measure is unique up to a multiplicative constant. Moreover, there exists a multiplier, $\Delta$, called the (*right*) *modular function* of the group, with the property that $\nu_r(dg) := \Delta(g^{-1})\nu_l(dg)$ is a right-Haar measure, which satisfies the obvious analogue of (10). Groups for which $\Delta(g) \equiv 1$; that is, for which right- and left-Haar measure are equivalent, are called *unimodular*. We now state two useful formulas that will be used repeatedly in the sequel (see Fremlin [4], Theorem 442K). If $\tilde{g} \in G$ and $h: G \to \mathbb{R}$ is an integrable function, then

(11) $$\int_G h(g\tilde{g}^{-1})\nu_l(dg) = \Delta(\tilde{g})\int_G h(g)\nu_l(dg)$$

and

$$\int_G h(g^{-1})\nu_l(dg) = \int_G h(g)\Delta(g^{-1})\nu_l(dg). \tag{12}$$

Now assume that $m(y) := \int_G f_Y(gy)j(g,y)\nu_l(dg)$ is positive for all $y \in \mathsf{Y}$ and finite for $\mu_y$-almost all $y \in \mathsf{Y}$. As in Section 4.2, let $N$ denote the $\mu_y$-null set of $y$ values for which $m(y) = \infty$ and set $\overline{\mathsf{Y}} = \mathsf{Y} \setminus N$. A routine calculation shows that, for $y \in \overline{\mathsf{Y}}$,

$$m(gy) = j(g^{-1}, y)\Delta(g^{-1})m(y). \tag{13}$$

This formula is basically equation (A1) from [10]. One consequence of (13) is that $gy \in \overline{\mathsf{Y}}$ for all $y \in \overline{\mathsf{Y}}$ and all $g \in G$. Let $Q$ be an operator on $L_0^2(f_Y)$ defined by

$$(Qh)(y) = \int_G \frac{h(gy)f_Y(gy)j(g,y)}{m(y)}\nu_l(dg)$$

when $y \in \overline{\mathsf{Y}}$ and $(Qh)(y) = \int_\mathsf{Y} h(y)f_Y(y)\mu_y(dy)$ when $y \in N$. This is the operator associated with the Markov chain on $\mathsf{Y}$ that evolves as follows. If the current state, $y$, is in $\overline{Y}$, then the distribution of the next state is that of $gy$ where $g$ is a random element from $G$ whose density (with respect to $\nu_l$) is $f_Y(gy)j(g,y)/m(y)$, and if $y \in N$, then the next state is from $f_Y$. Denote the chain and its Mtf by $\Psi = \{\Psi_n\}_{n=0}^\infty$ and $Q(y, dy')$.

PROPOSITION 3. *Suppose that $m(y)$ is positive for all $y \in \mathsf{Y}$ and finite for $\mu_y$-almost all $y \in \mathsf{Y}$ so that $Q$ is well-defined. Then the Mtf $Q$ is reversible with respect to $f_Y$.*

PROOF. As in the proof of Proposition 2, let $h_1, h_2 \in L_0^2(f_Y)$ be bounded. Then

$$\begin{aligned}\langle h_1, Qh_2\rangle &= \int_G \int_{\overline{\mathsf{Y}}} \frac{h_1(y)h_2(gy)f_Y(y)f_Y(gy)j(g,y)}{m(y)}\mu_y(dy)\nu_l(dg)\\
&= \int_G \left[\int_{\overline{\mathsf{Y}}} \frac{h_1(g^{-1}y)h_2(y)f_Y(g^{-1}y)f_Y(y)}{m(g^{-1}y)}\mu_y(dy)\right]\nu_l(dg)\\
&= \int_{\overline{\mathsf{Y}}} \frac{h_2(y)f_Y(y)}{m(y)}\left[\int_G h_1(g^{-1}y)f_Y(g^{-1}y)j(g^{-1},y)\Delta(g^{-1})\nu_l(dg)\right]\\
&\quad \times \mu_y(dy)\\
&= \int_{\overline{\mathsf{Y}}} \frac{h_2(y)f_Y(y)}{m(y)}\left[\int_G h_1(gy)f_Y(gy)j(g,y)\nu_l(dg)\right]\mu_y(dy)\\
&= \langle Qh_1, h_2\rangle,\end{aligned}$$

where the second through fourth equalities are due to, respectively, (8), (13) and (12). □

Compared with Theorem 1 in [10], our Proposition 3 is more general and has a stronger conclusion (reversibility versus invariance).

EXAMPLE 1 (continued). As noted previously, the multiplicative group is unimodular and $\nu_l(dg) = dg/g$ where $dg$ denotes Lebesgue measure on $\mathbb{R}^+$. Now, $m(y) = \int_G f_Y(gy)\chi(g)\nu_l(dg) = (2|y|)^{-1}$. Therefore, $N = \{0\}$ and $\overline{\mathsf{Y}}$ is the real line less the origin. For $y \neq 0$, $g \sim \text{Exp}(|y|)$ and for measurable $A \subset \overline{\mathsf{Y}}$, $Q(y, A) = \int_A q(y'|y)\mu_y(dy')$ where

$$q(y'|y) = e^{-|y'|}[I_{\mathbb{R}^+}(y)I_{\mathbb{R}^+}(y') + I_{\mathbb{R}^-}(y)I_{\mathbb{R}^-}(y')].$$

Again, the chain is not irreducible. However, for any fixed starting value in $\overline{\mathsf{Y}}$, the random variables $\Psi_1, \Psi_2, \Psi_3, \ldots$ are independent and identically distributed (i.i.d.). Indeed, if $\psi_0 > 0$, then $\Psi_1, \Psi_2, \Psi_3, \ldots$ are i.i.d. Exp(1) and if $\psi_0 < 0$ then $\Psi_1, \Psi_2, \Psi_3, \ldots$ are i.i.d. with common distribution equal to that of $-Z$ where $Z \sim \text{Exp}(1)$.

The behavior exhibited by $\Psi$ in the example above is not exceptional. Indeed, $Q$ has the special property that, conditional on any fixed starting value in $\overline{\mathsf{Y}}$, $\{\Psi_n\}_{n=1}^\infty$ is an i.i.d. sequence (which must be from $f_Y$ if the chain satisfies the usual regularity conditions). We will not prove this result here (due to space limitations), but we will prove that $Q$ is idempotent. For $n \in \mathbb{N} := \{1, 2, 3, \ldots\}$, let $Q^n(y, dy')$ denote the $n$-step Mtf.

PROPOSITION 4. *Suppose that $m(y)$ is positive for all $y \in \mathsf{Y}$ and finite for $\mu_y$-almost all $y \in \mathsf{Y}$ so that $Q$ is well-defined. For each $y \in \mathsf{Y}$, $Q^2(y, dy') = Q(y, dy')$ and hence $Q^n(y, dy') = Q(y, dy')$ for all $n \in \mathbb{N}$.*

PROOF. We prove the result in the case where $N = \varnothing$ and leave the extension to the general case to the reader. We will show that for $h \in L_0^2(f_Y)$, $(Q^2 h)(y) = (Q(Qh))(y) = (Qh)(y)$ for all $y \in \mathsf{Y}$. Indeed,

$$(Q(Qh))(y) = \int_G \left[\int_G \frac{h(g'gy)f_Y(g'gy)j(g',gy)}{m(gy)}\nu_l(dg')\right]\frac{f_Y(gy)j(g,y)}{m(y)}\nu_l(dg)$$

$$= \int_G \frac{f_Y(gy)}{m(y)m(gy)}\left[\int_G h(g'gy)f_Y(g'gy)j(g'g,y)\nu_l(dg')\right]\nu_l(dg)$$

$$= \int_G \frac{f_Y(gy)}{m(y)m(gy)}\left[\int_G \Delta(g^{-1})h(g'y)f_Y(g'y)j(g',y)\nu_l(dg')\right]\nu_l(dg)$$

$$= \int_G \int_G \frac{j(g,y)h(g'y)f_Y(gy)f_Y(g'y)j(g',y)}{m(y)m(y)}\nu_l(dg')\nu_l(dg)$$

$$= \int_G \frac{h(g'y)f_Y(g'y)j(g',y)}{m(y)m(y)} \left[ \int_G f_Y(gy)j(g,y)\nu_l(dg) \right] \nu_l(dg')$$
$$= (Qh)(y),$$

where the third and fourth equalities are due to, respectively, (11) and (13). □

The discussion preceding the statement of Proposition 4 suggests that it might be possible to use $Q$ to make i.i.d. draws from $f_Y$. Unfortunately, as we now explain, it is typically impossible to simulate $Q$ when the corresponding Markov chain is irreducible. Fix $y \in \mathsf{Y}$ and define

$$O_y = \{y' \in \mathsf{Y} : y' = gy \text{ for some } g \in G\}.$$

The set $O_y$ is called the *orbit* of $y$. The orbits induce an equivalence relation on the space $\mathsf{Y}$; that is, two points are equivalent if they are in the same orbit. Hence, $\mathsf{Y}$ can be *partitioned* into a collection of orbits. Clearly, when the Markov chain driven by $Q$ (or $Q_r$ for that matter) is started at the fixed value $y \in \overline{\mathsf{Y}}$, it remains forever in $O_y$. Therefore, if the probability measure associated with $f_Y$ puts positive mass on the complement of $O_y$, the Markov chain will not be $f_Y$-irreducible. Of course, the complement of $O_y$ definitely has measure zero when $O_y = \mathsf{Y}$; that is, when there is only one orbit. Unfortunately, when $\mathsf{Y}$ and $G$ are Euclidean spaces, the situations where there is only one orbit are those in which $g$ and $y$ have the same dimension. In practice, sampling from $f_Y(y)$ is not feasible and hence, if $g$ and $y$ share the same dimension, making draws from a density (in $g$) that is proportional to $f_Y(gy)j(g,y)$ will also likely be impossible. Loosely speaking, we are able to simulate $Q$ only when the corresponding Markov chain is reducible. While such reducible chains are not particularly useful by themselves, they can be used as part of a hybrid chain that is irreducible (see, e.g., Liu and Sabatti [10]) and they can be used to improve other chains such as the DA algorithm.

**5. General versions of PX-DA and Haar PX-DA.** Our general PX-DA algorithm has Mtd given by

$$p_r(x|x') = \int_\mathsf{Y} \int_\mathsf{Y} f_{X|Y}(x|y') Q_r(y, dy') f_{Y|X}(y|x') \mu_y(dy).$$

We now prove that $p_r$ is better than $p$ defined at (5) in both the efficiency ordering and the operator norm sense. We accomplish this by showing that $p_r$ is a DA algorithm. Let $P$ and $P_r$ denote the operators corresponding to $p$ and $p_r$.

PROPOSITION 5. *Let $r$ be a probability measure on $G$ such that $Q_r$ is well-defined. Then the Mtd $p_r$ is a DA algorithm. Thus, if $p$ and $p_r$ satisfy the usual regularity conditions, then $p_r \succeq_E p$ and $\|P_r\| \le \|P\|$.*

PROOF. Define $\tilde{f}(x,y) = \int_G f(x,gy) j(g,y) r(dg)$ and note that $\int_Y \tilde{f}(x,y) \times \mu_y(dy) = f_X(x)$. Hence, $\tilde{f}$ is a joint density on $X \times Y$ (with respect to $\mu_x \times \mu_y$) whose $x$ marginal is $f_X$. For $y \in \overline{Y}$,

$$\tilde{f}_{X|Y}(x|y) = \frac{\tilde{f}(x,y)}{\int_X \tilde{f}(x,y) \mu_x(dx)} = \frac{\int_G f(x,gy) j(g,y) r(dg)}{m_r(y)},$$

where, as in Section 4.2, $m_r(y) = \int_G f_Y(gy) j(g,y) r(dg)$. Also,

$$\tilde{f}_{Y|X}(y|x) = \frac{\tilde{f}(x,y)}{\int_Y \tilde{f}(x,y) \mu_y(dy)} = \frac{\int_G f(x,gy) j(g,y) r(dg)}{f_X(x)}$$

$$= \int_G f_{Y|X}(gy|x) j(g,y) r(dg).$$

Now,

$$\int_Y \tilde{f}_{X|Y}(x|y) \tilde{f}_{Y|X}(y|x') \mu_y(dy)$$

$$= \int_{\overline{Y}} \left[ \frac{1}{m_r(y)} \int_G f(x, g'y) j(g', y) r(dg') \right] \left[ \int_G f_{Y|X}(gy|x') j(g, y) r(dg) \right]$$
$$\times \mu_y(dy)$$

$$= \int_G \int_G \left[ \int_{\overline{Y}} \frac{1}{m_r(g^{-1}y)} f(x, g'g^{-1}y) j(g', g^{-1}y) f_{Y|X}(y|x') \mu_y(dy) \right]$$
$$\times r(dg') r(dg)$$

$$= \int_{\overline{Y}} \left[ \int_G \int_G \frac{f_{X|Y}(x|g'g^{-1}y) f_Y(g'g^{-1}y) j(g', g^{-1}y)}{m_r(g^{-1}y)} r(dg') r(dg) \right]$$
$$\times f_{Y|X}(y|x') \mu_y(dy)$$

$$= \int_Y \left[ \int_Y f_{X|Y}(x|y') Q_r(y, dy') \right] f_{Y|X}(y|x') \mu_y(dy) = p_r(x|x'),$$

where the second equality is due to (8) and the penultimate equality follows from the definition of $Q_r$. We conclude that $p_r$ is a DA algorithm. An appeal to Corollary 1 yields the result. □

In Proposition 5, the efficiency ordering result is new, but a special case of the operator norm result (where X, Y & $G$ are Euclidean spaces) is known— see L&W's Theorem 2.

Our general Haar PX-DA algorithm has Mtd given by

$$p^*(x|x') = \int_Y \int_Y f_{X|Y}(x|y')Q(y,dy')f_{Y|X}(y|x')\mu_y(dy),$$

where $Q(y,dy')$ is the Mtf defined in Section 4.3. Let $P^*$ denote the corresponding operator. Our next result establishes that the Haar PX-DA algorithm is better than *every* PX-DA algorithm in both the efficiency ordering and the operator norm sense. Before we state and prove the result, we explain the main idea. The most direct route to a proof would be to show that $Q \succeq_1 Q_r$ for every $r(dg)$, and then apply Theorem 3. However, we have not been able to establish that $Q \succeq_1 Q_r$. Alternatively, the reason we found success in comparing $p_r$ and $p$ is that $p_r$ is an improvement of the DA algorithm. At first glance, there is no such connection between $p^*$ and $p_r$. However, Proposition 5 says that $p_r$ is a DA algorithm and it turns out that $p^*$ can be represented as an improvement of $p_r$.

THEOREM 4. *Let $r$ be any probability measure on $G$ such that $Q_r$ is well-defined. Suppose that $m(y)$ is positive for all $y \in Y$ and finite for $\mu_y$-almost all $y \in Y$ so that $Q$ is well-defined. If $p_r$ and $p^*$ satisfy the usual regularity conditions, then $p^* \succeq_E p_r$ and $\|P^*\| \leq \|P_r\|$.*

PROOF. We prove the result in the case where $N = \varnothing$ (for both $m_r$ and $m$) and leave the extension to the general case to the reader. We know from Proposition 5 that $p_r$ is a DA algorithm with respect to the joint density $\tilde{f}(x,y) = \int_G f(x,gy)j(g,y)r(dg)$ and that $\int_X \tilde{f}(x,y)\mu_x(dx) = m_r(y)$. Let $\tilde{Q}$ be the Mtf on Y with invariant density $m_r(y)$ that is constructed according to the recipe in Section 4.3; that is, $\tilde{Q}$ is what we would have ended up with had we used $m_r(y)$ in place of $f_Y(y)$ in Section 4.3. We will show that

(14) $$p^*(x|x') = \int_Y \int_Y \tilde{f}_{X|Y}(x|y')\tilde{Q}(y,dy')\tilde{f}_{Y|X}(y|x')\mu_y(dy);$$

that is, $p^*$ is an improvement of $p_r$. First, if we substitute $m_r(y)$ for $f_Y(y)$ in the definition of $m(y)$, we have

$$\int_G m_r(gy)j(g,y)\nu_l(dg) = \int_G \left[\int_G f_Y(g'gy)j(g',gy)r(dg')\right]j(g,y)\nu_l(dg)$$

$$= \int_G \left[\int_G f_Y(g'gy)j(g'g,y)\nu_l(dg)\right]r(dg')$$

$$= \int_G f_Y(gy)j(g,y)\nu_l(dg) = m(y).$$

Hence, the function $m(y)$ is the same whether we use $f_Y$ or $m_r$. Now, using the definition of $\tilde{Q}$ and the calculation above, we have

$$\int_Y \tilde{f}_{X|Y}(x|y')\tilde{Q}(y,dy') = \frac{1}{m(y)}\int_G \tilde{f}_{X|Y}(x|g''y)m_r(g''y)j(g'',y)\nu_l(dg'').$$

Thus,

$$\int_{\mathsf{Y}} \bigg[\int_{\mathsf{Y}} \tilde{f}_{X|Y}(x|y')\tilde{Q}(y,dy')\bigg] \tilde{f}_{Y|X}(y|x')\mu_y(dy)$$

$$= \int_{\mathsf{Y}} \bigg[\int_G \bigg[\int_G \frac{f(x,g'g''y)j(g',g''y)r(dg')}{m_r(g''y)}\bigg] \frac{m_r(g''y)j(g'',y)\nu_l(dg'')}{m(y)}\bigg]$$

$$\times \bigg[\int_G f_{Y|X}(gy|x')j(g,y)r(dg)\bigg]\mu_y(dy)$$

$$= \int_G \int_G \int_G \bigg[\int_{\mathsf{Y}} \frac{f(x,g'g''y)j(g'g'',g^{-1}gy)j(g,y)f_{Y|X}(gy|x')}{m(y)}\mu_y(dy)\bigg]$$

$$\times \nu_l(dg'')r(dg')r(dg)$$

$$= \int_G \int_G \int_G \bigg[\int_{\mathsf{Y}} \frac{f(x,g'g''g^{-1}y)j(g'g'',g^{-1}y)f_{Y|X}(y|x')}{m(g^{-1}y)}\mu_y(dy)\bigg]$$

$$\times \nu_l(dg'')r(dg')r(dg)$$

$$= \int_G \int_G \int_{\mathsf{Y}} \bigg[\int_G \frac{f(x,g'g''g^{-1}y)j(g'g''g^{-1},y)\Delta(g^{-1})f_{Y|X}(y|x')}{m(y)}\nu_l(dg'')\bigg]$$

$$\times \mu_y(dy)r(dg')r(dg)$$

$$= \int_G \int_G \int_{\mathsf{Y}} \bigg[\int_G \frac{f(x,g'g''y)j(g'g'',y)f_{Y|X}(y|x')}{m(y)}\nu_l(dg'')\bigg]$$

$$\times \mu_y(dy)r(dg')r(dg)$$

$$= \int_G \int_G \int_{\mathsf{Y}} \bigg[\int_G \frac{f(x,g''y)j(g'',y)f_{Y|X}(y|x')}{m(y)}\nu_l(dg'')\bigg]\mu_y(dy)r(dg')r(dg)$$

$$= \int_{\mathsf{Y}} \bigg[\int_G \frac{f(x,g''y)j(g'',y)f_{Y|X}(y|x')}{m(y)}\nu_l(dg'')\bigg]\mu_y(dy)$$

$$= \int_{\mathsf{Y}} \bigg[\int_G \frac{f_{X|Y}(x|g''y)f_Y(g''y)j(g'',y)}{m(y)}\nu_l(dg'')\bigg]f_{Y|X}(y|x')\mu_y(dy)$$

$$= \int_{\mathsf{Y}} \bigg[\int_{\mathsf{Y}} f_{X|Y}(x|y')Q(y,dy')\bigg]f_{Y|X}(y|x')\mu_y(dy) = p^*(x|x'),$$

where the second equality follows from the properties of $j$, the third is from (8), the fourth is due to Fubini and (13), the fifth is a consequence of (11), the sixth is due to the left-invariance of $\nu_l$, the seventh follows from the fact that $r$ is a probability measure, and the penultimate equality is due to the definition of $Q$. Proposition 4 implies that $\tilde{Q}(y,dy')$ is idempotent and it follows from Proposition 1 that $p^*(x|x')$ is a DA algorithm. An application of Corollary 1 yields the result. $\square$

L&W proved that $\|P^*\| \le \|P_r\|$ in the special case where X, Y and $G$ are Euclidean spaces and $G$ is a unimodular group. Their proof relies heavily on a further assumption regarding the group structure that we now describe. Recall that Y can be partitioned into a set of orbits. A *cross section* is basically a subset of Y that intersects each orbit exactly once (see, e.g., Wijsman [24]). L&W assume the existence of a cross-section and a corresponding diffeomorphism that allows one to express each point in Y in terms of two quantities—its orbit and its position within its orbit. As L&W point out, the existence of a cross-section and an associated diffeomorphism is not guaranteed in general.

Recall from the discussion in Section 1 that L&W and M&vD developed (what appear to be) different strategies for handling the case in which $r$ is improper. We now demonstrate that the general Haar PX-DA Markov chain can be viewed as a marginal Markov chain associated with a nonpositive recurrent chain on a larger space. This result implies that, when the group structure is present, M&vD's chain (with left-Haar measure for the working prior) is exactly the same as L&W's Haar PX-DA algorithm. Suppose, as in most of the interesting applications, that $\nu_l(G) = \infty$. Following the ideas in M&vD, consider the function mapping $X \times Y \times G$ into $[0, \infty)$ that is defined by $\hat{f}(x, y, g) = f(x, gy)j(g, y)$. Now since

$$\int_G \int_Y \int_X \hat{f}(x,y,g)\mu_x(dx)\mu_y(dy)\nu_l(dg) = \nu_l(G),$$

$\hat{f}(x, y, g)$ is *not* integrable and therefore cannot be normalized to be a probability density function with respect to $\mu_x \times \mu_y \times \nu_l$. On the other hand, we can formally define "conditional" densities based on $\hat{f}$ as follows:

$$\hat{f}(y|x,g) = \frac{\hat{f}(x,y,g)}{\int_Y \hat{f}(x,y,g)\mu_y(dy)}$$

$$= \frac{f(x,gy)j(g,y)}{f_X(x)} = f_{Y|X}(gy|x)j(g,y),$$

and, for $y \in \overline{Y}$,

$$\hat{f}(x,g|y) = \frac{\hat{f}(x,y,g)}{\int_G \int_X \hat{f}(x,y,g)\mu_x(dx)\nu_l(dg)} = \frac{f(x,gy)j(g,y)}{m(y)}.$$

Therefore, despite the fact that $\hat{f}$ is not a density,

$$p^*((x,g)|(x',g')) = \int_Y \hat{f}(x,g|y)\hat{f}(y|x',g')\mu_y(dy)$$

is still a "DA-type" Mtd on $X \times G$. A routine calculation reveals that $f_X(x) \times \mu_x(dx)\nu_l(dg)$ is an invariant measure for the corresponding Markov chain, which we denote by $\{(X_n, G_n)\}_{n=0}^{\infty}$. However, $\int_G \int_\mathsf{X} f_X(x)\mu_x(dx)\nu_l(dg) = \nu_l(G)$, and hence the chain cannot be positive recurrent (Hobert [6]). On the other hand, the density of $X_{n+1}$ given $(X_n, G_n) = (x', g')$ is

$$\int_G p^*((x,g)|(x',g'))\nu_l(dg)$$
$$= \int_G \left[\int_\mathsf{\overline{Y}} \frac{f(x,gy)j(g,y)}{m(y)} f_{Y|X}(g'y|x')j(g',y)\mu_y(dy)\right]\nu_l(dg)$$
$$= \int_G \left[\int_\mathsf{\overline{Y}} \frac{f(x,gg'^{-1}y)j(g,g'^{-1}y)}{m(g'^{-1}y)} f_{Y|X}(y|x')\mu_y(dy)\right]\nu_l(dg)$$
$$= \int_\mathsf{\overline{Y}} \left[\int_G \frac{f(x,gg'^{-1}y)j(g,g'^{-1}y)}{m(y)j(g',y)\Delta(g')} f_{Y|X}(y|x')\nu_l(dg)\right]\mu_y(dy)$$
$$= \int_\mathsf{\overline{Y}} \left[\int_G \frac{f(x,gg'^{-1}y)j(gg'^{-1},y)\Delta(g'^{-1})}{m(y)} f_{Y|X}(y|x')\nu_l(dg)\right]\mu_y(dy)$$
$$= \int_\mathsf{\overline{Y}} \left[\int_G \frac{f(x,gy)j(g,y)}{m(y)} f_{Y|X}(y|x')\nu_l(dg)\right]\mu_y(dy) = p^*(x|x'),$$

where the second equality follows from (8), the third is due to Fubini and (13), the fourth is a consequence of the properties of $j$ and the fifth equality is due to (11). Since $\int_G p^*((x,g)|(x',g'))\nu_l(dg)$ does not depend on $g'$, it follows that $\{X_n\}_{n=0}^{\infty}$ itself is a Markov chain and the previous calculation shows that it is precisely the Markov chain driven by $p^*(x|x')$.

**Acknowledgments.** The authors thank Hani Doss, Galin Jones, Brett Presnell, Jeffrey Rosenthal, Vivekananda Roy and four reviewers for helpful comments and suggestions.

DEPARTMENT OF STATISTICS  
UNIVERSITY OF FLORIDA  
GAINESVILLE, FLORIDA 32611  
USA  
E-MAIL: jhobert@stat.ufl.edu

DEPARTMENT OF STATISTICS & CIS  
BARUCH COLLEGE, CUNY  
NEW YORK, NEW YORK 10010  
USA  
E-MAIL: dobrin_marchev@baruch.cuny.edu